\documentclass{amsart}

 \theoremstyle{definition}
 
 \theoremstyle{remark}
 
 \numberwithin{equation}{subsection}

\def\a{\alpha}
\def\b{\beta}

\def\D{\Delta}
\def\UU{{\mathcal U}}

\def\g{\gamma}

\def\Der{{\rm Der}}
\def\Inn{{\rm Inn}}

\def\Ker{{\rm Ker}}

\def\Im{{\rm Im}}

\def\v{\varphi}

\def\ssc{\scriptscriptstyle}

\def\cl{\centerline}

\def\rar{\rightarrow}

\def\vs{\vspace*}

\def\ni{\noindent}

\def\C{\mathbb{C}{\ssc\,}}
\def\F{\C}
\def\QED{\hfill$\Box$}
\numberwithin{equation}{section}
\newtheorem{theo}{Theorem}[section]
\newtheorem{defi}[theo]{Definition}

\newtheorem{lemm}[theo]{Lemma}
\newtheorem{prop}[theo]{Proposition}
\newtheorem{clai}{Claim}

\def\adddot{$\!\!\!${\bf.}\ \ }

\begin{document}


\cl {{\Large\bf Lie Bialgebra Structures on Generalized }}
 \cl {{\Large\bf
 \vs{10pt}  Heisenberg-Virasoro Algebra}\footnote{Supported by NSF of China (No.11001046), the Fundamental
Research Funds for the Central Universities, "Outstanding young
teachers of Donghua University" foundation.}}
\vs{6pt}

\cl{Hai Bo Chen$^{\,1)}$,  Ran Shen$^{\,2)}$, Jian Gang Zhang$^{\,
3)}$}


 \cl{\small
$^{1)}$,$^{2)}$College of Science, Donghua University, Shanghai,
201620, China}

\cl{\small $^{3)}$Department of Mathematics, Shanghai Normal
University, Shanghai, 200234, China}

\cl{\small Email:rshen@dhu.edu.cn}

\thanks{This work was Supported by NSF of China (No.11001046), the Fundamental
Research Funds for the Central Universities, "Outstanding young
teachers of Donghua University" foundation.}




\dedicatory{}


\vs{6pt}
\noindent{{\bf Abstract}  \,\,\,\,In this paper, Lie bialgebra structures on generalized
Heisenberg-Virasoro algebra $\mathfrak{L}$ are considered. Also, $H^1({\mathfrak{L}}
,\mathfrak{L}\bigotimes\mathfrak{L})$ is given explicitly. Moreover, it is
proved that all Lie bialgebra structure on centerless generalized
Heisenberg-Virasoro algebra $\overline{\mathfrak{L}}$ are coboundary
triangular.

\vs{6pt}
}

\vs{6pt}
\noindent{\bf MR(2000) Subject Classification}\,\,\,\, 17B62, 17B05, 17B37, 17B66


\section{Introduction}
\setcounter{section}{1}\setcounter{equation}{0} Lie bialgebra
structures on some Lie (super)algebras including generalized Witt
type, generalized Virasoro like type and generalized Weyl type Lie
algebras, the Schr$\ddot{o}$dinger-Virasoro Lie algebra, the $N=2$
superconformal algebra, etc., were constructed (cf.{\cite
{LCZ}},{\cite {HLS}},{\cite {M}},{\cite {SS}},{\cite {WSS1}},{\cite
{WSS2}},{\cite {YS1}},{\cite {YS2}}{\cite {YSS}}) since the notion
was first introduced by Drinfeld in 1983 (cf.{\cite {D1}},{\cite
{D2}}) in a connection with quantum groups. Recently, a general
method to obtain Lie (super)bialgebra structures on some Lie (super)
algebras related to Virasoro algebra was given in {\cite {LPZ}}.

In this paper, We will study the Lie bialgebra structures on Lie
algebra of generalized Heisenberg-Virasoro algebra which has been
studied in {\cite {LZ}}. It is the natural generalization of the
twisted Heisenberg-Virasoro algebra. The structure and
representations for generalized Heisenberg-Virasoro algebra were
studied in {\cite {LZ}},{\cite {SJS}}. However, Lie bialgebra
structures on generalized Heisenberg-Virasoro algebra have not yet
been considered.

Let us recall some definitions related to Lie bialgebras. For a
vector space $\mathfrak{L}$ over a field $\mathbb{F}$ of
characteristic zero. we define the {\it twist map} $\tau$ of
$\mathfrak{L}\otimes \mathfrak{L} $ and the {\it cyclic map} $\xi$
of $\mathfrak{L}\otimes \mathfrak{L} \otimes \mathfrak{L} $
\vs{-6pt}by
\begin{equation}
\tau:x\otimes y\mapsto y \otimes x, \ \ \ \ \ \ \xi: x \otimes
y\otimes z\mapsto y \otimes z \otimes x \mbox{ \ \ \ for \ \
}x,y,z\in \mathfrak{L} \vs{-6pt}.
\end{equation}
The definitions of a Lie algebra and Lie coalgebra can be
reformulated as follows. Then a {\it Lie algebra} can be defined as
a pair $(\mathfrak{L},\v)$ consisting of a vector space
$\mathfrak{L} $ and a bilinear map $\v :\mathfrak{L} \otimes
\mathfrak{L} \rar \mathfrak{L} $ (called the {\it bracket} of
$\mathfrak{L} $) satisfying the following conditions\vs{-6pt}:
\begin{eqnarray}
\label{Lie-s-s} \!\!\!\!\!\!\!\!\!\!\!\!&&
\Ker(1-\tau) \subset \Ker\,\v \mbox{ \ (skew-symmetry),}\\
\label{Lie-j-i} \!\!\!\!\!\!\!\!\!\!\!\!&& \v \cdot (1 \otimes \v )
\cdot (1 + \xi +\xi^{2}) =0 : \ \mathfrak{L}  \otimes \mathfrak{L}
\otimes \mathfrak{L}\rar \mathfrak{L}\mbox{ \ (Jacobi identity),}
\end{eqnarray}
denote $1$ is the identity map of $\mathfrak{L} \otimes \mathfrak{L}
$. Dually, a {\it Lie coalgebra} is a pair $(\mathfrak{L} , \D)$
consisting of a vector space $\mathfrak{L}$ and a linear map $\D:
\mathfrak{L} \to \mathfrak{L}\otimes \mathfrak{L}$ (called the {\it
cobracket} of $\mathfrak{L} $) satisfying the following conditions:
\begin{eqnarray}
\label{cLie-s-s} \!\!\!\!\!\!\!\!\!\!\!\!&&
\Im\,\D \subset \Im(1- \tau) \mbox{ \ (anti-commutativity),}\\
\label{cLie-j-i} \!\!\!\!\!\!\!\!\!\!\!\!&& (1 + \xi +\xi^{2}) \cdot
(1 \otimes \D) \cdot \D =0:\ \mathfrak{L}  \to \mathfrak{L} \otimes
\mathfrak{L} \otimes \mathfrak{L}\mbox{ \ (Jacobi identity).}
\end{eqnarray}
\begin{defi}
\label{def1} \rm A {\it Lie bialgebra} is a triple $(\mathfrak{L}
,\v, \D )$ such that $(\mathfrak{L} ,\v)$ is a Lie algebra and
$(\mathfrak{L} ,\D)$ is a Lie coalgebra and the following {\it
compatibility condition} holds:
\begin{eqnarray}
\label{bLie-d} \!\!\!\!\!\!\!\!\!\!\!\!&& \mbox{$\D  \v (x, y) = x
\cdot \D y - y \cdot \D x$ \ \ \ for \ \ $x, y \in \mathfrak{L} $, }
\end{eqnarray}
We shall use the symbol ``$\cdot$'' to stand for the {\it diagonal
adjoint action}:
\begin{equation}
\label{diag} x\cdot (\mbox{$\sum\limits_{i}$} {a_{i} \otimes b_{i}})
= \mbox{$\sum\limits_{i}$} ( {[x, a_{i}] \otimes b_{i} + a_{i}
\otimes [x, b_{i}]})
\end{equation}
for $x, a_{i}, b_{i} \in \mathfrak{L} $, and in general
$[x,y]=\v(x,y)$ for $x,y \in \mathfrak{L} $.
\end{defi}

\begin{defi}
\label{def2} \rm (1) A {\it coboundary Lie bialgebra} is a 4-tuple
$(\mathfrak{L} , \v, \D,r),$ where $(\mathfrak{L}, \v, \D)$ is a Lie
bialgebra and $r \in \Im(1 - \tau) \subset \mathfrak{L}  \otimes
\mathfrak{L} $ such that $\D=\D_r$ is a {\it coboundary of $r$},
where in general $\D_r$ (which is an inner derivation,
cf.~(\ref{inn})) is defined by,
\begin{equation}
\label{D-r} \D_r (x) = x \cdot r \mbox{ \ \ for \ \ }x \in
\mathfrak{L}.
\end{equation}

(2) A coboundary Lie bialgebra $(\mathfrak{L} , \v,\D, r)$ is {\it
triangular} if it satisfies the following {\it classical Yang-Baxter
Equation} (CYBE):
\begin{equation}
\label{CYBE} c(r)=0.
\end{equation}

(3) An element $r \in Im(1-\tau)
\subset\mathfrak{L}\otimes\mathfrak{L}$ is said to satisfy the {\it
modified Yang-Baxter Equation} (MYBE) if
\begin{equation}
\label{MYBE} x\cdot{c(r)}=0, \mbox{ \ \ for \ \ }x \in \mathfrak{L}.
\end{equation}
where $c(r)$ is defined by
\begin{equation}
\label{add1-} \mbox{$c(r) = [r^{12} , r^{13}] +[r^{12} , r^{23}]
+[r^{13} , r^{23}],$}
\end{equation}
and $r^{ij}$ are defined as follows: Denote $\UU(\mathfrak{L})$ the
universal enveloping algebra of $\mathfrak{L}$ and $1$ the identity
element of $\UU (\mathfrak{L})$. If $r =\sum_{i} {a_{i} \otimes
b_{i}} \in \mathfrak{L} \otimes \mathfrak{L} $, then $r^{ij}$ are
the following elements in $\UU (\mathfrak{L}) \otimes \UU
(\mathfrak{L} ) \otimes \UU(\mathfrak{L} )$:
$$\begin{array}{l}
 r^{12} = r\otimes 1=\sum \limits_{i}{a_{i} \otimes b_{i}
\otimes 1,} \\[12pt]
r^{13}= (1\otimes \tau)(\tau\otimes 1)=\sum \limits_{i} {a_{i}
\otimes 1 \otimes b_{i}}, \\[12pt]
r^{23} =1\otimes r= \sum \limits_{i}{1 \otimes a_{i} \otimes b_{i}}.
\end{array}$$
\end{defi}

The following results can be found in {\cite {D2}} and {\cite {NT}}.

\begin{lemm}
\label{some} {\rm(1)} For a Lie algebra ${\mathfrak{L}}$ and
$r\in\Im(1-\tau)\subset {\mathfrak{L}}$, the tripple
$({\mathfrak{L}},[\cdot,\cdot], \D)$ is a Lie bialgebra if and only
if $r$ satisfies MYBE.

{\rm(2)} For a Lie algebra ${\mathfrak{L}}$ and $r\in
\Im(1-\tau)\subset {\mathfrak{L}}$, we have
\begin{equation}
\label{add-c} (1 + \xi + \xi^{2}) \cdot (1 \otimes \D) \cdot \D (x)
= x \cdot c (r) \mbox{ \ for all \ }x\in \mathfrak{L}.
\end{equation}
\end{lemm}

Let us state our main results below. Suppose $\Gamma$ be  an abelian
group and $T$ is a vector space over $\mathbb{F}$. We always assume
that $T=\mathbb{F} \partial$ because of $dim T = 1$. The tensor
product $W=\mathbb{F} \Gamma \otimes_\mathbb{F} T$ is free left
$\mathbb{F}\Gamma$-module. We shall denote $t^x \partial=t^x \otimes
\partial$. Fix a pairing $\varphi: T \times \Gamma \rightarrow \mathbb{F}$,
which is $\mathbb{F}-$linear in the first variable and additive in
the second one. Denote: $\varphi (\partial, x) = \langle \partial, x
\rangle = \partial (x)$ for $x \in \Gamma$. If $\Gamma_0 = \{x\in
\Gamma: \partial (x)=0\}=0$ called $\varphi$ is nondegenerate. In
this paper,we always suppose that $\varphi$ is nondegenerate.

Clearly, from $\Gamma$ and $T$ is a vector space over $\mathbb{F}$,
the {\it generalized Heisenberg-Virasoro algebra}
$\mathfrak{L}:=\mathfrak{L}{(\Gamma)}$({\cite{LZ}}) is a Lie algebra
generated by $\{L_{x}=t^x\partial,I_{x}=t^x,C_L,C_I,C_{LI},x\in
\Gamma\}$, subject to the following relations:
$$\begin{array}{l}
[L_{x},L_{y}] = (y-x) L_{x+y}+ \delta_{x+y,0} \frac{1}{12}
(x^3 - x)C_L,\\[12pt]

[I_{x},I_{y}] = y \delta_{x+y,0} C_I,\\[12pt]

[L_{x},I_{y}]= y I_{x+y} + \delta_{x+y,0} (x^2- x)C_L,\\[12pt]

[\mathfrak{L},C_L] = [\mathfrak{L},C_I] = [\mathfrak{L},C_{LI}] = 0.

\end{array}$$
The Lie algebra $\mathfrak{L}$ has a generalized Heisenberg
subalgebra and a generalized Virasoro subalgebra interwined with a
2-cocycle. Set $\mathfrak{L}_x = \text{Span}_\mathbb{F}
\{L_{x},I_{x}\}$ for $x\in \Gamma \setminus\{0\}$, $\mathfrak{L}_0
=\text{Span}_\mathbb{F} \{L_{0},I_{0},C_L,C_I,C_{LI}\}$. Then
$\mathfrak{L}=\underset{x\in \Gamma}{\oplus} \mathfrak{L}_x $ is a
graded Lie algebra. Denote $\mathcal{C}$ the center of
$\mathfrak{L}$, then
$\mathcal{C}=\text{Span}_\mathbb{F}\{I_0,C_L,C_I,C_{LI}\}$.
  It is well known that the first cohomology group of
$\mathfrak{L}$ with coefficients in the module $V$ is isomorphic to
\begin{equation}
H^1(\mathfrak{L},V)\cong\Der(\mathfrak{L},V)/\Inn(\mathfrak{L},V),
\end{equation}
where $\Der(\mathfrak{L},V)$ is the set of \textit{derivations}
$d:\mathfrak{L}\to V$ which are linear maps satisfying
\begin{equation}
\label{deriv} d([x,y])=x\cdot d(y)-y\cdot d(x)\mbox{ \ for \ }x,y\in
\mathfrak{L},
\end{equation}
and $\Inn(\mathfrak{L},V)$ is the set of {\it inner derivations}
$a_{\rm inn},\,a\in V$, defined by
\begin{equation}
\label{inn} a_{\rm inn}:x\mapsto x\cdot a\mbox{ \ for \ }x\in
\mathfrak{L}.
\end{equation}

\section{Lie bialgebra structures on the generlized
Heisenberg-Virasoro algebra}

\setcounter{section}{2} \setcounter{theo}{0}\setcounter{equation}{0}

\begin{defi}
\label{def3} For any  $\lambda \in \mathbb{F}$ ,$C\in \mathcal {C}$,
we define the map $\lambda \otimes C:\mathfrak{L}\rightarrow
\mathfrak{L}\otimes\mathfrak{L}$ by
\begin{eqnarray}\label{equa2.8}
\begin{split}
 (\lambda \otimes
C)(\omega_\alpha)= & \lambda(1-\delta_{\alpha,0})w_1 I_\alpha
\otimes C,&
\end{split}
\end{eqnarray}
$\mbox{ \ for  \ }
\omega_\alpha=w_1L_\alpha+w_2I_\alpha+\delta_{\alpha,0} Z \in
 \mathfrak{L}_\alpha,
 \mbox{ \ where  \ } Z\in Span_\mathbb{F}\{C_L,C_{LI},C_I\},
 \alpha\in \Gamma.$
\end{defi}

Obviously, for any $C\in \mathcal {C}$, $\lambda \in \mathbb{F}$,
$\lambda \otimes C \in\Der(\mathfrak{L},\mathfrak{L}\otimes
\mathfrak{L})$ and it is an outer derivation. Furthermore,
$\mathbb{F}\otimes C=\{\lambda\otimes C, \lambda\in \mathbb{F}\}$ is
a subalgebra of $\Der(\mathfrak{L},\mathfrak{L}\otimes
\mathfrak{L})$, denoted by $\mathbb{F}\otimes \mathbb{F} C$.
Similarly, we have the derivation subalgebra $\mathbb{F}C\otimes
\mathbb{F}$.

The main result of this section is
\begin{theo}
\label{main} {\rm(1)} Let $(\mathfrak{L},[\cdot,\cdot],\Delta)$ be a
Lie bialgebra such that $\Delta$ has the decomposition
$\Delta_r+\sigma$  with respect to $\Der({\mathfrak{L}} ,V)=
\Inn({\mathfrak{L}} ,V)\oplus (\mathbb{F}\otimes
\mathbb{F}C+\mathbb{F}C\otimes \mathbb{F})$, where $r\in
V(\mathrm{mod}  \mathcal {C}\otimes \mathcal {C})$ and $\sigma\in
\mathbb{F} \otimes \mathbb{F}C+\mathbb{F}C\otimes \mathbb{F},
\sigma(\mathfrak{L})\subseteq \mathrm{Im}(1-\tau)$. Then, $r\in
\mathrm{Im}(1-\tau)$. Furthermore,
$(\mathfrak{L},[\cdot,\cdot],\sigma)$ is a Lie bialgebra.

 {\rm(2)}
An element $r \in \Im(1 - \tau) \subset {\mathfrak{L}}\otimes
{\mathfrak{L}} $
 satisfies CYBE in $(\ref{CYBE})$ if and only if it
satisfies MYBE in $(\ref{MYBE})$.

{\rm(3)} Regarding $V={\mathfrak{L}} \otimes {\mathfrak{L}}$ as an
${\mathfrak{L}}$-module under the adjoint diagonal action of
${\mathfrak{L}} $ in $(\ref{diag})$, we have $H^1({\mathfrak{L}}
,V)=\Der({\mathfrak{L}} ,V)/\Inn({\mathfrak{L}} ,V)\cong
\mathbb{F}\otimes \mathbb{F}C+\mathbb{F}C\otimes \mathbb{F}$.
\end{theo}

We give the proof of Theorem \ref{main} by several lemmas and
propositions.

At first, Theorem \ref{main}(2) follows from the following result.

\begin{lemm}
\label{lemma2} Denote by ${\mathfrak{L}}^{\otimes n}$ the tensor
product of $n$ copies of ${\mathfrak{L}} $. Regarding
${\mathfrak{L}}^{\otimes n}$ as an ${\mathfrak{L}}$-module under the
adjoint diagonal action of ${\mathfrak{L}}$, suppose
$c\in{\mathfrak{L}}^{\otimes n}$ satisfying $a\cdot c=0$ for all
$a\in{\mathfrak{L}} $. Then $c=0$.
\end{lemm}
\ni{\it Proof.} The lemma is obtained by using the same arguments as
in the proof of Lemma 2.2 in {\cite{WSS2}}.

Theorem \ref{main}(3) follows from the following proposition.

\begin{prop}
\label{lemma3} $\Der({\mathfrak{L}} ,V)=\Inn({\mathfrak{L}}
,V)+\mathbb{F}\otimes \mathbb{F}C+\mathbb{F}C\otimes \mathbb{F}$,
where $V= {\mathfrak{L}} \otimes {\mathfrak{L}}$.
\end{prop}

\ni{\it Proof.} We shall divide the proof of the proposition into
several claims. Note that $V=\oplus_{\a\in \Gamma}V_\a$ is
$\Gamma$-graded with $V_\a=\sum_{\b+\g=\a}
{\mathfrak{L}}_\b\otimes{\mathfrak{L}}_\g$, where
${\mathfrak{L}}_\a=\F L_\a\oplus \mathbb{C}I_\alpha\oplus
\delta_{\alpha,0}(\mathbb{C}C_L+\mathbb{C}C_I+\mathbb{C}C_{LI})$ for
$\a\in \Gamma$. A derivation $D\in\Der({\mathfrak{L}} ,V)$ is {\it
homogeneous of degree $\a\in \Gamma$} if $D(V_\b) \subset V_{\a
+\b}$ for all $\b \in \Gamma$. Denote $\Der({\mathfrak{L}} , V)_\a =
\{D\in \Der({\mathfrak{L}} , V) \,|\,{\rm deg\,}D =\a\}$ for $\a\in
\Gamma$.

\begin{clai}
\label{clai1} \rm Every derivation $D\in\Der({\mathfrak{L}} ,V)$.
Then
\begin{equation}
\label{summable} \mbox{$D = \sum\limits_{\a \in \Gamma} D_\a ,
\mbox{  \ where \ }D_\a \in \Der({\mathfrak{L}}, V)_\a,$}
\end{equation}
which holds in the sense that for every $\omega \in {\mathfrak{L}}
$, only finitely many $D_\a (\omega)\neq 0,$ and $D(\omega) =
\sum_{\a \in \Gamma} D_\a(\omega)$ (we call such a sum in
(\ref{summable}) {\it summable}).
\end{clai}
\ \indent For $\a\in \Gamma$, we define a homogeneous linear map
$D_\a : {\mathfrak{L}}\rightarrow V$ of degree $\alpha$ as follows:
For any $\omega\in {\mathfrak{L}} _\b$ with $\b\in \Gamma$, write
$d(\omega)=\sum_{\g\in \Gamma}v_\g\in V$ with $v_\g\in V_\g$, then
we set $D_\a(\omega)=v_{\a+\b}$. Obviously
$D_\a\in\Der({\mathfrak{L}} ,V)_\a$ and we have (\ref{summable}).
\begin{clai}
\label{clai2} \rm If $0\neq\a\in \Gamma$, then
$D_\a\in\Inn({\mathfrak{L}} ,V)$.
\end{clai}
\ \indent  Denote $u=\alpha^{-1}D_\alpha(L_0)\in V_\alpha.$ For any
$\omega_\beta\in \mathfrak{L}_\beta, \beta\in \Gamma,$  applying
$D_\alpha$ to $[L_0,\omega_\beta]=\beta \omega_\beta,$  using
$D_\alpha(\omega_j)\in V_{\alpha+\beta}$ and the action of $L_0$ on
$V_{\alpha+\beta}$ is the scalar $\alpha+\beta$, we have
\begin{equation}\label{equa-add-1}
(\a+\b)D_{\a}(\omega_\beta) - \omega_\beta\cdot D_{\a}(L_0)=\b
D_{\a}(\omega_\beta),
\end{equation}
i.e., $D_{\a}(\omega_\beta)=u_{\rm
inn}(\omega_\beta)(\mathrm{cf}.(\ref{inn})).$ Thus $D_{\a}=u_{\rm
inn}$ is inner.

\begin{clai}
\label{clai3} \rm $D_0(L_0)=D_0(C)=0(\mbox{ mod }\mathbb{F}(\mathcal
{C}\otimes \mathcal {C}))$.
\end{clai}

In order to prove this, applying $D_0$ to $[L_0,\omega]=\b \omega$
for $\omega \in {\mathfrak{L}}_\beta, \b \in \Gamma$, as in
(\ref{equa-add-1}), we obtain that \mbox{$\omega\cdot
D_0(L_0)=0(\mbox{ mod }\mathbb{F}(\mathcal {C}\otimes \mathcal
{C}))$.} Thus by lemma \ref{lemma2}, $D_0(L_0)=0$. Similarly, by
applying $D_0$ to $[C,\omega]=0$, we obtain $D_0(C)=0 (\mbox{ mod
}\mathbb{F}(\mathcal {C}\otimes \mathcal {C}))$. \vskip4pt

\begin{clai}
\label{clai3} \rm By replacing $D_0$ by $D_0-u_{\rm
inn}-(\lambda\otimes C+C \otimes \eta)$ for some $u\in
V_0,\lambda,\eta\in \mathbb{F}$, we can suppose
$D_0(\mathfrak{L}_{\mu})\equiv0(\mbox{ mod }\mathbb{F}(\mathcal
{C}\otimes \mathcal {C})$ for $\mu \in\Gamma$.
\end{clai}

We can write,under modulo $\mathbb{F}(\mathcal {C}\otimes \mathcal
{C})$ (where $\mu\in \Gamma/\mathbb{Z}$ means $\mu$ is a
representative of the coset $\mu+\mathbb{Z}$ in $\Gamma$, in case
$\mu\in \mathbb{Z}$ or $\mu \in \Gamma/\mathbb{Z}$ we always choose
$\mu=0$),
\begin{eqnarray}\label{equa2.9}
\begin{split}
D_0(L_{\pm 1})=&\underset{\mu\in \Gamma/\mathbb{Z},i\in
\mathbb{Z}}{\sum}{a^{\pm}_{\mu,i} L_{\mu+i\pm 1} \otimes
L_{-\mu-i}}+  \underset{\mu\in \Gamma/\mathbb{Z},i\in
\mathbb{Z}}{\sum}{b^{\pm}_{\mu,i}
L_{\mu+i\pm 1} \otimes I_{-\mu-i}}\\
&+\underset{\mu\in \Gamma/\mathbb{Z},i\in
\mathbb{Z}}{\sum}{c^{\pm}_{\mu,i} I_{\mu+i\pm 1} \otimes
L_{-\mu-i}}+  \underset{\mu\in \Gamma/\mathbb{Z},i\in
\mathbb{Z}}{\sum}{d^{\pm}_{\mu,i} I_{\mu+i\pm 1} \otimes
I_{-\mu-i}}\\ &+ a^{\pm }L_{\pm1}\otimes C+b^{\pm }C \otimes
L_{\pm1} +c^{\pm }I_{\pm1}\otimes C +d^{\pm }C \otimes I_{\pm1} ,
\end{split}
\end{eqnarray}
 and we set  $b^{+}_{0,0}=c^{+}_{0,-1}=d^{+}_{0,0}=
 d^{+}_{0,-1}=b^{-}_{0,0}=c^{-}_{0,1}=
 d^{-}_{0,0}=d^{-}_{0,1}=0$, for some
 $ a^{\pm}_{\mu,i},b^{\pm}_{\mu,i},c^{\pm}_{\mu,i},d^{\pm}_{\mu,i},
  a^{\pm},  \\
 b^{\pm },c^{\pm },d^{\pm } \in
  \mathbb{F} $, where $\{(\mu,i)\mid a^{\pm }_{\mu,i}\neq 0 \}$,  $\{(\mu,i)\mid b^{\pm }_{\mu,i}\neq 0
  \}$,$\{(\mu,i)\mid c^{\pm }_{\mu,i}\neq 0 \}$ and $\{(\mu,i)\mid d^{\pm}_{\mu,i}\neq 0 \}$
  are finite sets.In the following,to simplify notations, we shall
  always omit the superscript $``^+"$; for example,
  $a_{\mu,i}=a^+_{\mu,i}$. Note that for $\mu \in \Gamma/\mathbb{Z},i\in
  \mathbb{Z}$, we have
  $$\begin{array}{l}
  (L_{\mu+i} \otimes L_{-\mu-i})_{inn}(L_1)=(\mu+i-1)L_{\mu+i+1}\otimes L_{-\mu-i}
  -(\mu+i+1)L_{\mu+i}\otimes L_{-\mu-i+1},  \\[12pt]
  (L_{\mu+i} \otimes I_{-\mu-i})_{inn}(L_1)=(\mu+i-1)L_{\mu+i+1}\otimes I_{-\mu-i}
  -(\mu+i)L_{\mu+i}\otimes I_{-\mu-i+1},  \\[12pt]
  (I_{\mu+i} \otimes L_{-\mu-i})_{inn}(L_1)=(\mu+i)I_{\mu+i+1}\otimes L_{-\mu-i}
  -(\mu+i+1)I_{\mu+i}\otimes L_{-\mu-i+1},  \\[12pt]
 (I_{\mu+i} \otimes I_{-\mu-i})_{inn}(L_1)=(\mu+i)I_{\mu+i+1}\otimes I_{-\mu-i}
 -(\mu+i)I_{\mu+i}\otimes I_{-\mu-i+1},  \\[12pt]
 (L_0 \otimes C)_{inn}(L_1)=-L_1\otimes C,  \\[12pt]
  (C \otimes L_0)_{inn}(L_1)=-C\otimes L_1.  \\[12pt]
 \end{array}$$
Denote
$$\begin{array}{l}
 M_\mu=max\{|i|\mid a_{\mu,i}\neq 0 \},
 N_\mu=max\{|i|\mid b_{\mu,i}\neq 0 \}, \\[12pt]
E_\mu=max\{|i|\mid c_{\mu,i}\neq 0 \},
F_\mu=max\{|i|\mid d_{\mu,i}\neq 0\}.  \\[12pt]
\end{array}$$
 Using the above equations and the induction on
 $M_\mu+N_\mu+E_\mu+F_\mu$, by replacing $D_0$ by $D_0-u_{inn}$,
 where $u$ is a combination of some $L_{\mu+i} \otimes L_{-\mu-i}$, $L_{\mu+i} \otimes I_{-\mu-i}$,
 $I_{\mu+i} \otimes L_{-\mu-i}$, $I_{\mu+i} \otimes I_{-\mu-i}$, we
 can suppose
  $$\begin{array}{l}
  a= b=0,   \\[12pt]
  a_{\mu,i} =0  \mbox{ if } \mu=0; i\neq -2, 1,  \mbox{ or } \mu\neq 0,i\neq 0,\\[12pt]
  b_{\mu,i} =0  \mbox{ if } \mu=0; i\neq -1, 1,  \mbox{ or } \mu\neq 0,i\neq 0,\\[12pt]
  c_{\mu,i} =0  \mbox{ if } \mu=0; i\neq -2, 0,  \mbox{ or } \mu\neq 0,i\neq 0,\\[12pt]
  d_{\mu,i} =0  \mbox{ if } \mu=0; i\neq -1, 0,  \mbox{ or } \mu\neq 0,i\neq 0.\\[12pt]
  \end{array}$$
Thus we have, under modulo $\mathbb{F}(\mathcal {C}\otimes \mathcal
{C})$,

\begin{eqnarray}\label{equa2.10}
\begin{split}
 D_0(L_1)=& a_{0,-2}L_{-1}\otimes L_2+ a_{0,1}L_2\otimes L_{-1}+
 \underset{0\neq\mu\in \Gamma/\mathbb{Z}}{\sum}{a_{\mu,i} L_{\mu+1} \otimes
 L_{-\mu}}  \\&
 + b_{0,-1}L_0\otimes I_1 + b_{0,1}L_2\otimes I_{-1}
 +\underset{0\neq\mu\in \Gamma/\mathbb{Z}}{\sum}{b_{\mu,0}
 L_{\mu+1} \otimes I_{-\mu}}    \\&
 + c_{0,-2}I_{-1}\otimes L_2 + c_{0,0}I_1\otimes
 L_0 +\underset{0\neq\mu\in \Gamma/\mathbb{Z}}{\sum}{c_{\mu,0}I_{\mu+1} \otimes L_{-\mu}}
 \\&
 +d_{0,-1} I_0 \otimes I_1+d_{0,0} I_1 \otimes I_0+\underset{0\neq\mu\in \Gamma/\mathbb{Z}}{\sum}{d_{\mu,0}I_{\mu+1} \otimes I_{-\mu}}\\&
+c I_{1}\otimes C +d C \otimes I_{1}.
\end{split}
\end{eqnarray}

Applying $D_0$  to $[L_{-1},L_1]=2L_0$, we obtain, under modulo
$\mathbb{F} (\mathcal {C}\otimes \mathcal {C})$,
$$\begin{array}{ll}
3a_{0,-2}L_{-1}\otimes L_1+ 3a_{0,1}L_1\otimes L_{-1}+
 \underset{0\neq\mu\in \Gamma/\mathbb{Z}}{\sum}{a_{\mu,0}[(\mu+2)L_\mu \otimes
 L_{-\mu}}-(\mu-1) L_{\mu+1}\otimes L_{-\mu-1}]   \\[12pt]
  + b_{0,-1}[L_{-1}\otimes I_1+L_0\otimes I_0+2L_0\otimes C_L] + b_{0,1}[3L_1\otimes
  I_{-1} -L_2\otimes I_{-2}]     \\[12pt]
  +\underset{0\neq\mu\in \Gamma/\mathbb{Z}}{\sum}{b_{\mu,0}
[(\mu+2) L_\mu \otimes I_{-\mu}} {-\mu L_{\mu+1}\otimes I_{-\mu-1}]}
 + c_{0,-2}[-I_{-2}\otimes L_2+3I_{-1}\otimes L_1 ]   \\[12pt]
 + c_{0,0}[I_0 \otimes
 L_0+2C_L\otimes L_0+I_1\otimes L_{-1}]
 +\underset{0\neq\mu\in \Gamma/\mathbb{Z}}{\sum}c_{\mu,0}[(\mu+1)I_\mu
  \otimes L_{-\mu}-(\mu-1)I_{\mu+1}\otimes L_{-\mu-1}]   \\[12pt]
+\underset{0\neq\mu\in
\Gamma/\mathbb{Z}}{\sum}{d_{\mu,0}[(\mu+1)I_\mu \otimes
I_{-\mu}-\mu I_{\mu+1}\otimes I_{-\mu-1}}]     \\[12pt]
=\underset{\mu\in \Gamma/\mathbb{Z},i\in
\mathbb{Z}}{\sum}{a^{-}_{\mu,i} [(\mu+i-2)L_{\mu+i} \otimes L_{-\mu-i}-(\mu+i+1)L_{\mu+i-1}\otimes L_{-\mu-i+1}]}  \\[12pt]
+\underset{\mu\in \Gamma/\mathbb{Z},i\in
\mathbb{Z}}{\sum}{b^{-}_{\mu,i} [(\mu+i-2)L_{\mu+i} \otimes I_{-\mu-i}-(\mu+i)L_{\mu+i-1}\otimes I_{-\mu-i+1}]}   \\[12pt]
+\underset{\mu\in \Gamma/\mathbb{Z},i\in
 \mathbb{Z}}{\sum}{c^{-}_{\mu,i} [(\mu+i-1)I_{\mu+i} \otimes L_{-\mu-i}-(\mu+i+1)I_{\mu+i-1}\otimes L_{-\mu-i+1}]}   \\[12pt]
+ \underset{\mu\in \Gamma/\mathbb{Z},i\in
 \mathbb{Z}}{\sum}{d^{-}_{\mu,i} [(\mu+i-1)I_{\mu+i} \otimes I_{-\mu-i}-(\mu+i)I_{\mu+i-1}\otimes I_{-\mu-i+1}]}    \\[12pt]
-2a^{-} L_0\otimes C-2 b^{-} C\otimes L_0.
\end{array}$$

Comparing the coefficients of $L_{\mu+i} \otimes L_{-\mu-i}$  with
$\mu\in \Gamma/\mathbb{Z},i\in \mathbb{Z}$, we obtain
$$\begin{array}{lll}
(\mu+i-2)a^{-}_{\mu,i}-(\mu+i+2)a^{-}_{\mu,i+1}    \\[12pt]
=3\delta_{\mu,0}\delta_{i,-1}a_{0,-2}+3\delta_{\mu,0}\delta_{i,1}a_{0,1}
+(1-\delta_{\mu,0})\delta_{i,0}(\mu+2)a_{\mu,0}+(1-\delta_{\mu,0})\delta_{i,1}(1-\mu)a_{\mu,0}.
\end{array}$$
Since $\{(\mu,i)\mid a^{-}_{\mu,i}\neq 0 \}$  is a finite set, fix
$\mu\in \Gamma/\mathbb{Z},i\in \mathbb{Z}$,     we obtain
$$\begin{array}{llll}
a^{-}_{0,m+1}=a^{-}_{0,-m}=a^{-}_{\mu,i}=0 \mbox{ \ for \ } m\geq 2,
\mu\neq 0, i\neq 1,
\end{array}$$
and we have the following relations:
\begin{equation}
\label{equa2.11} a^{-}_{0,-1}=-a_{0,-2}-\frac{1}{3}a^{-}_{0,0},
a^{-}_{0,1}= -a^{-}_{0,0},
a^{-}_{0,2}=-a_{0,1}+\frac{1}{3}a^{-}_{0,0}, a^{-}_{\mu,1}
=-a_{\mu,0} \mbox{ \ for \ } \mu\neq 0.
\end{equation}
Comparing the coefficients of $L_{\mu+i} \otimes I_{-\mu-i}$  with
$\mu\in \Gamma/\mathbb{Z},i\in \mathbb{Z}$, we obtain
$$\begin{array}{llll}
 (\mu+i-2)b^{-}_{\mu,i}-(\mu+i+1)b^{-}_{\mu,i+1}    \\[12pt]
 =\delta_{\mu,0}(\delta_{i,-1}+\delta_{i,0})b_{0,-1}+\delta_{\mu,0}(3\delta_{i,1}-\delta_{i,2})b_{0,1}
 +(1-\delta_{\mu,0})\delta_{i,0}(\mu+2)b_{\mu,0}-(1-\delta_{\mu,0})\delta_{i,1}\mu
 b_{\mu,0}.
 \end{array}$$
Since $\{(\mu,i)\mid b^{-}_{\mu,i}\neq 0 \}$  is a finite set, fix
$\mu\in \Gamma/\mathbb{Z},i\in \mathbb{Z}$,and $b^{-}_{0,0}=0$, we
obtain
 $$\begin{array}{llll}
 b_{0,-1}=b_{0,1}=b_{\mu,0}
 =b^{-}_{\mu,i}=b^{-}_{0,i}=0  \mbox{ \ for \ } \mu\neq 0.
 \end{array}$$
Comparing the coefficients of $I_{\mu+i} \otimes L_{-\mu-i}$  with
  $\mu\in \Gamma/\mathbb{Z},i\in \mathbb{Z}$, we obtain
  $$\begin{array}{llll}
  (\mu+i-1)c^{-}_{\mu,i}-(\mu+i+2)c^{-}_{\mu,i+1}    \\[12pt]
  =\delta_{\mu,0}(-\delta_{i,-2}+3\delta_{i,-1})c_{0,-2}+\delta_{\mu,0}(\delta_{i,1}+\delta_{i,0})c_{0,0}
  +(1-\delta_{\mu,0})\delta_{i,0}(\mu+1)c_{\mu,0}+(1-\delta_{\mu,0})\delta_{i,1}(1-\mu)
  c_{\mu,0}.
  \end{array}$$
  Since $\{(\mu,i)\mid c^{-}_{\mu,i}\neq 0 \}$  is a finite set,
  fix $\mu\in \Gamma/\mathbb{Z},i\in \mathbb{Z}$, and $c^{-}_{0,1}=0$,   we obtain
  $$\begin{array}{llll}
  c_{0,-2}=c_{0,0}=c_{\mu,0}
  =c^{-}_{\mu,i}=c^{-}_{0,i}=0  \mbox{ \ for \ } \mu\neq 0.
  \end{array}$$
Comparing the coefficients of $I_{\mu+i} \otimes I_{-\mu-i}$  with
  $\mu\in \Gamma/\mathbb{Z},i\in \mathbb{Z}$, we obtain
  $$\begin{array}{llll}
  (\mu+i-1)d^{-}_{\mu,i}-(\mu+i+1)c^{-}_{\mu,i+1}    \\[12pt]
  =(1-\delta_{\mu,0})\delta_{i,0}(\mu+1)d_{\mu,0}-(1-\delta_{\mu,0})\delta_{i,1}\mu
  d_{\mu,0}.
  \end{array}$$
   Since $\{(\mu,i)\mid d^{-}_{\mu,i}\neq 0 \}$  is a finite set,
   fix $\mu\in \Gamma/\mathbb{Z},i\in \mathbb{Z}$, and $d^{-}_{0,0}=d^{-}_{0,-1}=0$,   we obtain
   $$\begin{array}{llll}
    d^{-}_{\mu,m}=d^{-}_{0,i}=0  \mbox{ \ for \ } m\neq 1,\mu\neq 0,
   \end{array}$$
 and we have the following relation:
 \begin{equation}
 \label{equa2.12}
 d^{-}_{\mu,1}= -d_{\mu,0}  \mbox{ \ for \ } \mu\neq 0.
 \end{equation}
Comparing the coefficients of $L_0 \otimes C,C\otimes L_0$, we
obtain
 $$\begin{array}{lllll}
 a^{-}= b^{-}=0.
 \end{array}$$
Consequently, under modulo $\mathbb{F}(\mathcal {C}\otimes \mathcal
{C})$, we can rewrite
\begin{eqnarray}\label{equa2.13}
\begin{split}
 D_0(L_{-1})=& a^{-}_{0,-1}L_{-2}\otimes L_1+ a^{-}_{0,0}L_{-1}\otimes L_0
 + a^{-}_{0,1}L_0\otimes L_{-1}+ a^{-}_{0,2}L_1\otimes
 L_{-2} \\&
 -\underset{0\neq\mu\in \Gamma/\mathbb{Z}}{\sum}{a_{\mu,0} L_\mu \otimes
 L_{-\mu-1}}-\underset{0\neq\mu\in \Gamma/\mathbb{Z}}{\sum}{d_{\mu,0}I_\mu \otimes I_{-\mu-1}}
+c^{-}I_{-1}\otimes C \\& +d^{-}C \otimes I_{-1},
\end{split}
\end{eqnarray}
where the coefficients satisfy (\ref{equa2.11}) and
(\ref{equa2.12}).

Under modulo $\mathbb{F}(\mathcal {C}\otimes \mathcal {C})$, we can
write
 \begin{eqnarray}\label{equa2.14}
 \begin{split}
 D_0(L_{\pm 2})=&\underset{\mu\in \Gamma/\mathbb{Z},i\in
 \mathbb{Z}}{\sum}{a^{'\pm}_{\mu,i} L_{\mu+i\pm 2} \otimes
 L_{-\mu-i}}+  \underset{\mu\in \Gamma/\mathbb{Z},i\in
 \mathbb{Z}}{\sum}{b^{'\pm}_{\mu,i}
 L_{\mu+i\pm 2} \otimes I_{-\mu-i}}\\
 &+\underset{\mu\in \Gamma/\mathbb{Z},i\in
 \mathbb{Z}}{\sum}{c^{'\pm}_{\mu,i} I_{\mu+i\pm 2} \otimes
 L_{-\mu-i}}+  \underset{\mu\in \Gamma/\mathbb{Z},i\in
 \mathbb{Z}}{\sum}{d^{'\pm}_{\mu,i} I_{\mu+i\pm 2} \otimes
 I_{-\mu-i}}\\ &+ a^{'\pm}L_{\pm2}\otimes
 C+b^{'\pm}C \otimes L_{\pm2}
 +c^{'\pm}I_{\pm2}\otimes C  +d^{'\pm}C \otimes
 I_{\pm2},
 \end{split}
\end{eqnarray}
and we set $b^{'+}_{0,0}=c^{'+}_{0,-2}=d^{'+}_{0,0}=
d^{'+}_{0,-2}=b^{'-}_{0,0}=c^{'-}_{0,2}=
d^{'-}_{0,0}=d^{'-}_{0,2}=0$, for some
$a^{'\pm}_{\mu,i},b^{'\pm}_{\mu,i},c^{'\pm}_{\mu,i},d^{'\pm}_{\mu,i},
 a^{'\pm}, \\
 b^{'\pm},c^{'\pm},d^{'\pm}\in
 \mathbb{F}$, where $\{(\mu,i)\mid a^{'\pm}_{\mu,i}\neq 0 \}$,  $\{(\mu,i)\mid b^{'\pm}_{\mu,i}\neq 0
 \}$,$\{(\mu,i)\mid c^{'\pm}_{\mu,i}\neq 0 \}$ and $\{(\mu,i)\mid d^{'\pm}_{\mu,i}\neq 0 \}$
 are finite sets.Now we shall omit the superscript $``^+ "$ again; for example,
  $a^{'}_{\mu,i}=a^{'+}_{\mu,i}$.

Applying $D_0$  to $[L_{-1},L_2]=3L_1$, we have, under modulo
$\mathbb{F} (\mathcal {C}\otimes \mathcal {C})$,
$$\begin{array}{llll}
&\underset{\mu\in \Gamma/\mathbb{Z},i\in
\mathbb{Z}}{\sum}{a^{'}_{\mu,i} [(\mu+i+3)L_{\mu+i+ 1} \otimes
L_{-\mu-i}-(\mu+i-1)L_{\mu+i+2}\otimes L_{-\mu-i-1}]}   \\&

+\underset{\mu\in \Gamma/\mathbb{Z},i\in
\mathbb{Z}}{\sum}{b^{'}_{\mu,i} [(\mu+i+3)L_{\mu+i+1} \otimes
I_{-\mu-i}-(\mu+i)L_{\mu+i+2}\otimes I_{-\mu-i-1}]}    \\&

+\underset{\mu\in \Gamma/\mathbb{Z},i\in
\mathbb{Z}}{\sum}{c^{'}_{\mu,i} [(\mu+i+2)I_{\mu+i+1} \otimes
L_{-\mu-i}-(\mu+i-1)I_{\mu+i+2}\otimes L_{-\mu-i-1}]}     \\&
+\underset{\mu\in \Gamma/\mathbb{Z},i\in
\mathbb{Z}}{\sum}{d^{'}_{\mu,i} [(\mu+i+2)I_{\mu+i+1} \otimes
I_{-\mu-i}-(\mu+i)I_{\mu+i+2}\otimes I_{-\mu-i-1}]}     \\&
+3a^{'}L_1\otimes C+3b^{'}C \otimes L_{1} +2c^{'}I_1\otimes C
+2d^{'}C \otimes I_1 +a^{-}_{0,-1}[4L_0\otimes L_1+L_{-2}\otimes
L_3]   \\&

+ a^{-}_{0,0} [3L_1\otimes L_0+2L_{-1}\otimes L_2]
 + a^{-}_{0,1}[2L_2\otimes L_{-1}+3L_0\otimes L_1]
 + a^{-}_{0,2}[L_3\otimes L_{-2}+4L_1\otimes L_0]       \\&
 +\underset{0\neq\mu\in \Gamma/\mathbb{Z}}{\sum}{a^{-}_{\mu,0} [(\mu-2)L_{\mu+2} \otimes
 L_{-\mu-1}-(\mu+3)L_\mu \otimes L_{-\mu-3}]}          \\&

 +\underset{0\neq\mu\in \Gamma/\mathbb{Z}}{\sum}
 {d^{-}_{\mu,0}[\mu I_{\mu+2} \otimes I_{-\mu-1}-(\mu+1)I_\mu \otimes I_{-\mu+1}]}
+c^{-}I_{1}\otimes C+d^{-}C \otimes I_{1} \\&

=3a_{0,-2}L_{-1}\otimes L_2+ 3a_{0,1}L_2\otimes L_{-1}+
\underset{0\neq\mu\in \Gamma/\mathbb{Z}}{\sum}{3a_{\mu,0} L_{\mu+1}
\otimes L_{-\mu}}  \\&

+\underset{0\neq\mu\in \Gamma/\mathbb{Z}}{\sum}{3d_{\mu,0}I_{\mu+1}
\otimes I_{-\mu}} +3cI_{1}\otimes C +3d C\otimes I_{1}.
\end{array}$$

Comparing the coefficients of $L_{\mu+i+1} \otimes L_{-\mu-i}$  with
$\mu\in \Gamma/\mathbb{Z},i\in \mathbb{Z}$, we obtain
$$\begin{array}{llll}
3\delta_{i,-2}a_{0,-2}+3\delta_{i,1}a_{0,1}  =& (i+3)a^{'}_{0,i}
-(i-2)a^{'}_{0,i-1}
+\delta_{i,-3}a^{-}_{0,-1}+2\delta_{i,-2}a^{-}_{0,0}
 \\&  +
\delta_{i,-1}(4a^{-}_{0,-1} +3a^{-}_{0,1}) +
 \delta_{i,0}(4a^{-}_{0,2} +3a^{(L_{-1})}_{0,0}) +
 2\delta_{i,1}a^{-}_{0,1}+\delta_{i,2}a^{-}_{0,2},
\end{array}$$

and

$$\begin{array}{llll}
 (\mu+i+3)a^{'}_{\mu,i}-(\mu+i-2)a^{'}_{\mu,i-1}
+\delta_{i,1}(\mu-2)a_{\mu,0}-\delta_{i,-1}(\mu+3)a_{\mu,0}
=3\delta_{i,0}a_{\mu,0} \mbox{ \ for \ } \mu\neq 0.
\end{array}$$

Using (\ref{equa2.11}), we obtain

$$\begin{array}{llll}
a^{-}_{0,-1}=a^{-}_{0,2}=a^{'}_{0,m}=a^{'}_{0,-m-2}=
a_{\mu,0}=a^{'}_{\mu,i}=0 \mbox{ \ for \ } m\geq 2, \mu\neq
0,i\in\mathbb{Z},
\end{array}$$

and

 \begin{eqnarray}\label{equa2.15}
 \begin{split}
  &a^{-}_{0,0}=-a^{-}_{0,1}=3a_{0,1}= -3a_{0,-2}, \\&
 a^{'}_{0,-3}= a^{'}_{0,1}-6a_{0,1},
 a^{'}_{0,-2}=-4a^{'}_{0,1}+15a_{0,1}, \\&
 a^{'}_{0,-1}= 6a^{'}_{0,1}-18a_{0,1},
 a^{'}_{0,0}= -4a^{'}_{0,1}+9a_{0,1}.
 \end{split}
\end{eqnarray}

Comparing the coefficients of $L_{\mu+i+1} \otimes I_{-\mu-i}$  with
$\mu\in \Gamma/\mathbb{Z},i\in \mathbb{Z}$, we obtain

$$\begin{array}{llll}
b^{'}_{\mu,i}=0.
\end{array}$$

Comparing the coefficients of $I_{\mu+i+1} \otimes L_{-\mu-i}$  with
$\mu\in \Gamma/\mathbb{Z},i\in \mathbb{Z}$, we obtain

$$\begin{array}{llll}
c^{'}_{\mu,i}=0.
\end{array}$$

Comparing the coefficients of $I_{\mu+i+1} \otimes I_{-\mu-i}$  with
$\mu\in \Gamma/\mathbb{Z},i\in \mathbb{Z}$, we obtain

$$\begin{array}{llll}
(i+2)d^{'}_{0,i} -(i-1)d^{'}_{0,i-1}=0 ,
\end{array}$$

and

$$\begin{array}{llll}
 (\mu+i+2)d^{'}_{\mu,i}-(\mu+i-1)d^{'}_{\mu,i-1}
+\mu\delta_{i,1}d_{\mu,0}-(\mu+1)\delta_{i,-1}d_{\mu,0}
=3\delta_{i,0}d_{\mu,0} \mbox{ \ for \ } \mu\neq 0.
\end{array}$$

Using $d^{'}_{0,0}=d^{'}_{0,-2}=0$, we obtain

$$\begin{array}{llll}
d^{'}_{0,i}=d^{'}_{\mu,m}=d^{'}_{\mu,-m-1}=0 \mbox{ \ for \ } m\geq
1, \mu\neq 0,i\in\mathbb{Z},
\end{array}$$

and

 \begin{eqnarray}\label{equa2.16}
 \begin{split}
 d^{'}_{\mu,-1}=d^{'}_{\mu,0}=d_{\mu,0}, \mbox{ \ for \ } \mu\neq 0.
\end{split}
\end{eqnarray}

Now (\ref{equa2.10}) and (\ref{equa2.13}), respectively,become,under
modulo $\mathbb{F}(\mathcal {C}\otimes \mathcal {C}),$

\begin{eqnarray}\label{equa2.17}
 \begin{split}
D_0(L_{1})=& a_{0,1}(-L_{-1}\otimes L_2+L_2\otimes L_{-1})
\\& +\underset{0\neq\mu\in
\Gamma/\mathbb{Z}}{\sum}{d_{\mu,0}I_{\mu+1} \otimes I_{-\mu}} +c
I_{1}\otimes C +d C \otimes I_{1},
\end{split}
\end{eqnarray}
\begin{eqnarray}\label{equa2.18}
 \begin{split}
D_0(L_{-1})=& 3a_{0,1}(L_{-1}\otimes L_0-L_0\otimes L_{-1})
\\& -\underset{0\neq\mu\in
\Gamma/\mathbb{Z}}{\sum}{d_{\mu,0}I_\mu \otimes I_{-\mu-1}}
+c^{-}I_{-1}\otimes C +d^{-}C \otimes I_{-1}.
\end{split}
\end{eqnarray}

Applying $D_0$ to $[L_{-2},L_1]=3L_{-1},$ we have, under modulo
$\mathbb{F}(\mathcal {C}\otimes \mathcal {C}),$

$$\begin{array}{llll}
 & a_{0,1}[-L_{-3}\otimes L_2-4L_{-1}\otimes L_0-
 \frac{1}{2}L_{-1}\otimes C_L+4L_0\otimes L_{-1}+\frac{1}{2} C_L \otimes L_{-1}+L_2\otimes
 L_{-3}]  \\&

 +\underset{0\neq\mu\in
\Gamma/\mathbb{Z}}{\sum}{d_{\mu,0}[(\mu+1)I_{\mu-1} \otimes
I_{-\mu}}-\mu I_{\mu+1}\otimes I_{-\mu-2}] +c I_{-1}\otimes C +d C
\otimes I_{-1}  \\&

-(\underset{\mu\in \Gamma/\mathbb{Z},i\in
\mathbb{Z}}{\sum}{a^{'-}_{\mu,i} [(\mu+i-3)L_{\mu+i-1} \otimes
L_{-\mu-i}-(\mu+i+1)L_{\mu+i-2}\otimes L_{-\mu-i+1}]}   \\&

+\underset{\mu\in \Gamma/\mathbb{Z},i\in
\mathbb{Z}}{\sum}{b^{'-}_{\mu,i} [(\mu+i-3)L_{\mu+i-1} \otimes
I_{-\mu-i}-(\mu+i)L_{\mu+i-2}\otimes I_{-\mu-i+1}]}    \\&

+\underset{\mu\in \Gamma/\mathbb{Z},i\in
\mathbb{Z}}{\sum}{c^{'-}_{\mu,i} [(\mu+i-2)I_{\mu+i-1} \otimes
L_{-\mu-i}-(\mu+i+1)I_{\mu+i-2}\otimes L_{-\mu-i+1}]}     \\&
+\underset{\mu\in \Gamma/\mathbb{Z},i\in
\mathbb{Z}}{\sum}{d^{'-}_{\mu,i} [(\mu+i-2)I_{\mu+i-1} \otimes
I_{-\mu-i}-(\mu+i)I_{\mu+i-2}\otimes I_{-\mu-i+1}]}     \\& -3a^{'-}
L_{-1}\otimes C-3b^{'-} C \otimes L_{-1} -2c^{'-} I_{-1}\otimes
C-2d^{'-} C \otimes I_{-1})
\\&

=9a_{0,1}(L_{-1}\otimes L_0-L_0\otimes L_{-1})-3
 \underset{0\neq\mu\in \Gamma/\mathbb{Z}}{\sum}{d_{\mu,0}I_\mu \otimes
I_{-\mu-1}}

+3c^{-} I_{-1}\otimes C +3d^{-} C \otimes I_{-1}.
\end{array}$$

Comparing the coefficients, Using $d^{'-}_{0,0}=d^{'-}_{0,2}=0$,we
obtain
$$\begin{array}{llll}
a_{0,1}=a^{'-}_{\mu,i}=a^{'-}_{0,m+2}=a^{'-}_{0,-m}=a^{'-}=0 \mbox{
\ for \ } m\geq 2,\mu\neq 0, i\in\mathbb{Z},
\end{array}$$

$$\begin{array}{llll}
b^{'-}_{\mu,i}=c^{'-}_{\mu,i}=b^{'-}=d^{'-}_{0,i}=0 \mbox{ \ for \ }
\mu\neq 0. i\in\mathbb{Z},
\end{array}$$

$$\begin{array}{llll}
d^{'-}_{\mu,m}=d^{'-}_{\mu,-m+3}=0 \mbox{ \ for \ } m\geq 3, \mu\neq
0,
\end{array}$$

and

$$\begin{array}{llll}
a^{'-}_{0,0}=a^{'-}_{0,2}=-4a^{'-}_{0,-1},
a^{'-}_{0,1}=6a^{'-}_{0,-1}, a^{'-}_{0,3}=a^{'-}_{0,-1} ,
\end{array}$$

\begin{eqnarray}\label{equa3.18}
 \begin{split}
d^{'-}_{\mu,1}=d^{'-}_{\mu,2}=-d_{\mu,0} \mbox{ \ for \ } \mu \neq
0.
\end{split}
\end{eqnarray}

Now (\ref{equa2.17}),(\ref{equa2.18}), and (\ref{equa2.14}),
respectively, become, under modulo $\mathbb{F}(\mathcal {C} \otimes
\mathcal {C}),$
\begin{eqnarray}\label{equa2.19}
 \begin{split}
D_0(L_{1})=& \underset{0\neq\mu\in
\Gamma/\mathbb{Z}}{\sum}{d_{\mu,0}I_{\mu+1} \otimes I_{-\mu}} +c
I_{1}\otimes C +d C \otimes I_{1},
\end{split}
\end{eqnarray}

 \begin{eqnarray}\label{equa2.20}
 \begin{split}
D_0(L_{-1})= \ -\underset{0\neq\mu\in
\Gamma/\mathbb{Z}}{\sum}{d_{\mu,0}I_\mu \otimes I_{-\mu-1}}
+c^{-}I_{-1}\otimes C +d^{-}C \otimes I_{-1},
\end{split}
\end{eqnarray}

\begin{eqnarray}\label{equa2.21}
 \begin{split}
D_0(L_{2})=& a^{'}_{0,1}(L_{-1}\otimes L_3-4L_0\otimes
L_2+6L_{1}\otimes L_1-4L_2\otimes L_0+L_3\otimes L_{-1})  \\&
 +\underset{0\neq\mu\in
\Gamma/\mathbb{Z}}{\sum}{d_{\mu,0}(I_{\mu+2} \otimes
I_{-\mu}+I_{\mu+1} \otimes I_{-\mu+1})} +c^{'}I_{2}\otimes C +d^{'}C
\otimes I_{2},
\end{split}
\end{eqnarray}

\begin{eqnarray}\label{equa2.22}
 \begin{split}
D_0(L_{-2})=& a^{'-}_{0,-1}(L_{-3}\otimes L_1-4L_{-2}\otimes
L_0+6L_{-1}\otimes L_{-1}-4L_0\otimes L_{-2}+L_1\otimes L_{-3}) \\&
 -\underset{0\neq\mu\in
\Gamma/\mathbb{Z}}{\sum}{d_{\mu,0}(I_{\mu-1} \otimes
I_{-\mu-1}+I_{\mu} \otimes I_{-\mu-2})} +c^{'-}I_{-2}\otimes C
+d^{'-}C \otimes I_{-2}.
\end{split}
\end{eqnarray}
Note that
$$\begin{array}{llll}
(I_{\mu-1}\otimes I_{-\mu+1})_{inn}(L_2)=(\mu-1)I_{\mu+1}\otimes
I_{-\mu+1}-(\mu-1)I_{\mu-1}\otimes I_{-\mu+3},   \\

(I_{\mu-3}\otimes I_{-\mu+3})_{inn}(L_{2})=(\mu-3)I_{\mu-1}\otimes
I_{-\mu+3}-(\mu-3)I_{\mu-3}\otimes I_{-\mu+5}.
\end{array}$$
Using these two equations, by replacing $D_0$ by $D_0-u_{inn}$,
where $u$ is a combination of $I_{\mu-1}\otimes I_{-\mu+1}$ and
$I_{\mu-3}\otimes I_{-\mu+3}$ (this replacement does not affect the
above equations (\ref{equa2.19}),(\ref{equa2.20})and
(\ref{equa2.22}) ), under modulo $\mathbb{F}(\mathcal {C}\otimes
\mathcal {C})$, we can rewrite (\ref{equa2.21}) as
\begin{eqnarray}\label{equa3.21}
 \begin{split}
D_0(L_{2})=& a^{'}_{0,1}(L_{-1}\otimes L_3-4L_0\otimes
L_2+6L_{1}\otimes L_1-4L_2\otimes L_0+L_3\otimes L_{-1})  \\&
 +\underset{0\neq\mu\in
\Gamma/\mathbb{Z}}{\sum}{d_{\mu,0}(I_{\mu+2} \otimes
I_{-\mu}+I_{\mu-3} \otimes I_{-\mu+5})} +c^{'}I_{2}\otimes C +d^{'}C
\otimes I_{2},
\end{split}
\end{eqnarray}

Applying $D_0$ to $[L_{-2},L_2]=L_0+\frac{1}{2}C_L,$ we obtain
$$\begin{array}{llll}
a^{'}_{0,1}=-a^{'-}_{0,-1}, d_{\mu,0}=0  \mbox{ \ for \ } \mu\neq 0.
\end{array}$$
Using (\ref{equa2.12}),(\ref{equa2.16})and(\ref{equa3.18}), we
 obtain
$$\begin{array}{llll}
d^{-}_{\mu,1}=d^{'}_{\mu,-1}=d^{'}_{\mu,0}=d^{'-}_{\mu,1}
=d^{'-}_{\mu,2}=0 \mbox{ \ for \ } \mu\neq 0.
\end{array}$$
Denote
$$\begin{array}{llll}
u= L_{-1}\otimes L_1-2L_0\otimes L_0+L_1\otimes L_{-1}.
\end{array}$$
Replacing $D_0$ by $D_0+a^{'}_{0,1}u_{inn}$(this replacement does
not affect the above two  equations
(\ref{equa2.19})and(\ref{equa2.20})), we obtain
\begin{eqnarray}\label{equa2.23}
 \begin{split}
D_0(L_{2})=c^{'}I_{2}\otimes C +d^{'}C \otimes I_{2},
\end{split}
\end{eqnarray}

\begin{eqnarray}\label{equa2.24}
 \begin{split}
D_0(L_{-2})= c^{'-}I_{-2}\otimes C +d^{'-}C \otimes I_{-2}.
\end{split}
\end{eqnarray}

At last, we can rewrite (\ref{equa2.9}) and (\ref{equa2.14}) as
\begin{eqnarray}\label{equa2.25}
\begin{split}
D_0(L_{\pm 1})= c^{\pm}I_{\pm1}\otimes C +d^{\pm}C \otimes I_{\pm1}
,
\end{split}
\end{eqnarray}

\begin{eqnarray}\label{equa2.26}
 \begin{split}
 D_0(L_{\pm 2})=
 c^{'\pm}I_{\pm2}\otimes C  +d^{'\pm}C \otimes
 I_{\pm2}.
 \end{split}
\end{eqnarray}

Write under modulo $\mathbb{F}(\mathcal {C} \otimes \mathcal {C})$,
\begin{eqnarray}\label{equa2.27}
\begin{split}
D_0(I_{\pm 1})=&\underset{\mu\in \Gamma/\mathbb{Z},i\in
\mathbb{Z}}{\sum}{e^{\pm}_{\mu,i} L_{\mu+i\pm 1} \otimes
L_{-\mu-i}}+  \underset{\mu\in \Gamma/\mathbb{Z},i\in
\mathbb{Z}}{\sum}{f^{\pm}_{\mu,i}
L_{\mu+i\pm 1} \otimes I_{-\mu-i}}\\
&+\underset{\mu\in \Gamma/\mathbb{Z},i\in
\mathbb{Z}}{\sum}{g^{\pm}_{\mu,i} I_{\mu+i\pm 1} \otimes
L_{-\mu-i}}+  \underset{\mu\in \Gamma/\mathbb{Z},i\in
\mathbb{Z}}{\sum}{h^{\pm}_{\mu,i} I_{\mu+i\pm 1} \otimes
I_{-\mu-i}}\\ &+ e^{\pm}L_{\pm1}\otimes C+f^{\pm}C \otimes L_{\pm1}
+g^{\pm}I_{\pm1}\otimes C +h^{\pm}C \otimes I_{\pm1} ,
\end{split}
\end{eqnarray}
and we set $f^{+}_{0,0}=g^{+}_{0,-1}=h^{+}_{0,0}=
 h^{+}_{0,-1}=f^{-}_{0,0}=g^{-}_{0,1}=
 h^{-}_{0,0}=h^{-}_{0,1}=0$, for some
 $e^{\pm}_{\mu,i},f^{\pm}_{\mu,i},g^{\pm}_{\mu,i},h^{\pm}_{\mu,i},
 \\
  e^{\pm},f^{\pm},g^{\pm},h^{\pm}\in
  \mathbb{F}$, where $\{(\mu,i)\mid e^{\pm}_{\mu,i}\neq 0 \}$,  $\{(\mu,i)\mid f^{\pm}_{\mu,i}\neq 0
  \}$,$\{(\mu,i)\mid g^{\pm}_{\mu,i}\neq 0 \}$ and $\{(\mu,i)\mid h^{\pm}_{\mu,i}\neq 0 \}$
  are finite sets.In the following,to simplify notations, we shall
  always omit the superscript $``^+"$; for example,
  $e_{\mu,i}=e^{+}_{\mu,i}$. Note that for $\mu \in \Gamma/\mathbb{Z},i\in
  \mathbb{Z}$, we have

 $$\begin{array}{l}
(L_0 \otimes C)_{inn}(I_1)=-I_1\otimes C,  \\[12pt]
(C \otimes L_0)_{inn}(I_1)=-C\otimes I_1.  \\[12pt]
 \end{array}$$
Denote
      $M^{'}_\mu=max\{|i|\mid e_{\mu,i}\neq 0 \}$,  $N^{'}_\mu=max\{|i|\mid f_{\mu,i}\neq 0\}$,
      $E^{'}_\mu=max\{|i|\mid g_{\mu,i}\neq 0 \}$,  $F^{'}_\mu=max\{|i|\mid h_{\mu,i}\neq 0\}$.

Thus we have, under modulo $\mathbb{F}(\mathcal {C}\otimes \mathcal
{C})$,

\begin{eqnarray}\label{equa2.28}
\begin{split}
D_0(I_{1})=&\underset{\mu\in \Gamma/\mathbb{Z},i\in
\mathbb{Z}}{\sum}{e_{\mu,i} L_{\mu+i+ 1} \otimes L_{-\mu-i}}+
\underset{\mu\in \Gamma/\mathbb{Z},i\in \mathbb{Z}}{\sum}{f_{\mu,i}
L_{\mu+i+1} \otimes I_{-\mu-i}}\\
&+\underset{\mu\in \Gamma/\mathbb{Z},i\in
\mathbb{Z}}{\sum}{g_{\mu,i} I_{\mu+i+ 1} \otimes L_{-\mu-i}}+
\underset{\mu\in \Gamma/\mathbb{Z},i\in \mathbb{Z}}{\sum}{h_{\mu,i}
I_{\mu+i+ 1} \otimes I_{-\mu-i}}\\ &+ e L_{1}\otimes C+f C \otimes
L_{1}
\end{split}
\end{eqnarray}
Applying $D_0$ to $[L_{-1},I_1]=0$, we obtain
\begin{eqnarray}\label{equa2.29}
\begin{split}
D_0(I_{1})=0,
\end{split}
\end{eqnarray}
By the same method, we can obtain
\begin{eqnarray}\label{equa2.30}
\begin{split}
D_0(I_{-1})=0.
\end{split}
\end{eqnarray}
Using (\ref{equa2.25}),(\ref{equa2.26}),(\ref{equa2.29}) and
(\ref{equa2.30}), under modulo $\mathbb{F}(\mathcal {C}\otimes
\mathcal {C})$, we have
\begin{eqnarray}\label{equa2.31}
\begin{split}
D_0(L_{ i})= c_{i}I_{ i}\otimes C +d_{i}C \otimes I_{ i} \mbox { \
for \ } i\in \mathbb{Z},
\end{split}
\end{eqnarray}
\begin{eqnarray}\label{equa2.32}
\begin{split}
D_0(I_{i})=0  \mbox { \ for \ } i\in \mathbb{Z}.
\end{split}
\end{eqnarray}
For a fixed $0\neq \mu \in \Gamma/\mathbb{Z},$  Under modulo
$\mathbb{F}(\mathcal {C} \otimes \mathcal {C})$, we can write
 \begin{eqnarray}\label{equa2.33}
 \begin{split}
 D_0(L_{\mu})=&\underset{\nu\in \Gamma/\mathbb{Z},j\in
 \mathbb{Z}}{\sum}{a_{\nu,j} L_{\mu+\nu+j} \otimes
 L_{-\nu-j}}+\underset{\nu\in \Gamma/\mathbb{Z},j\in
 \mathbb{Z}}{\sum}{b_{\nu,j} L_{\mu+\nu+j} \otimes
 I_{-\nu-j}}   \\&
 +\underset{\nu\in \Gamma/\mathbb{Z},j\in
 \mathbb{Z}}{\sum}{c_{\nu,j} I_{\mu+\nu+j} \otimes
 L_{-\nu-j}}+\underset{\nu\in \Gamma/\mathbb{Z},j\in
 \mathbb{Z}}{\sum}{d_{\nu,j} I_{\mu+\nu+j} \otimes
 I_{-\nu-j}}   \\&
 + a_{\mu}L_{\mu}\otimes
 C+b_{\mu}C \otimes L_{\mu}
 +c_{\mu}I_{\mu}\otimes C  +d_{\mu}C \otimes
 I_{\mu},
 \end{split}
\end{eqnarray}
for some
 $a_{\mu,j},b_{\mu,j},c_{\mu,j},d_{\mu,j}  \in
  \mathbb{F}$, where $A^{'}=\{(\nu,j)\mid a_{\nu,j}\neq 0 \}$,  $B^{'}=\{(\nu,j)\mid b_{\nu,j}\neq 0
  \}$, $C^{'}=\{(\nu,j)\mid c_{\nu,j}\neq 0 \}$, $D^{'}=\{(\nu,j)\mid d_{\nu,j}\neq 0
  \}$ are finite sets. Let $i\gg 0$. Applying $D_0$ to
  $[L_i,[L_{-i},L_\mu]=(\mu+i)(\mu-2i)L_\mu$, we obtain
\begin{eqnarray}\label{equa2.34}
\begin{split}
L_i\cdot L_{-i}\cdot D_0(L_\mu)=(\mu+i)(\mu-2i)D_0(L_\mu)+L_i\cdot
L_\mu \cdot D_0(L_{-i})\\+(\mu+i)L_{\mu-i}\cdot D_0(L_i)
\end{split}
\end{eqnarray}
Define the total order on $\Gamma\times \Gamma$ by
$$\begin{array}{l}
(\alpha,\beta)>(\alpha^{'},\beta^{'}) \Leftrightarrow \alpha>
\alpha^{'} \mbox{ \ or \ } \alpha = \alpha^{'}, \beta> \beta^{'}.
\end{array}$$

If $A^{'}\neq \emptyset,$ let $(\nu_0,j_0)$ be the maximal element
in $A{'}$, in this case, $L_{\mu+\nu_0+j_0+i}\otimes
L_{-\nu_0-j_0-i}$ is the leading term of $L_i\cdot L_{-i}\cdot
D_0(L_\mu)$, a contradiction to equation (\ref{equa2.32}).
Similarly, we have $B^{'}= C^{'}=D^{'}=\emptyset$. Thus, we can
rewrite (\ref{equa2.33}) as under modulo $\mathbb{F}(\mathcal
{C}\otimes \mathcal{C})$,

$$\begin{array}{l}
 D_0(L_{\mu})=
 a_{\mu}L_{\mu}\otimes
 C+b_{\mu}C \otimes L_{\mu}
 +c_{\mu}I_{\mu}\otimes C  +d^{(L_{\mu})}_{\mu}C \otimes
 I_{\mu}.
 \end{array}$$
Then we can denote $u=\frac{a_\mu}{\mu}L_0\otimes
C+\frac{b_\mu}{\mu}C\otimes L_0$, replacing $D_0$ by $D_0+u_{inn}$,
we obtain
\begin{eqnarray}\label{equa2.35}
\begin{split}
 D_0(L_{\mu})=c_{\mu}I_{\mu}\otimes C  +d_{\mu}C \otimes
 I_{\mu}.
\end{split}
\end{eqnarray}
 Similarly, we have $
 D_0(I_{\nu})=e_{\nu}L_{\nu}\otimes C  +f_{\nu}C \otimes
 L_{\nu}.$

Applying $D_0$ to $[L_{\nu-1},I_1]=I_\nu(0\neq \nu \in
\Gamma/\mathbb{Z})$, we obtain, under modulo $\mathbb{F} (\mathcal
{C}\otimes \mathcal {C})$, $D_0(I_{\nu})=0.$ Using (\ref{equa2.31})
and (\ref{equa2.35}), we have $ D_0(L_{\mu})=c_{\mu}I_{\mu}\otimes C
+d_{\mu}C \otimes
 I_{\mu} \mbox{ \ for \ } \mu\in \Gamma\setminus \{0\}$, applying $D_0$ to
 $[L_{\mu},L_\nu]=(\nu-\mu)L_{\mu+\nu}$, we obtain $c_\mu,d_\mu \in \mathbb{F}$.
 Claim $4$ is proved.

\begin{clai}\adddot \label{clai4} \rm
We can suppose $D_0=0$ by replacing $D_0$ by
$D_0-u_{inn}-(\lambda\otimes C+C \otimes \eta)$ for some $u\in V_0,
\lambda, \eta\in \mathbb{F}.$
\end{clai}
\vskip4pt
\par

Replacing $D_0$ by $D_0-u_{inn}-(\lambda\otimes C+C \otimes \eta)$
for some $u\in V_0,\lambda, \eta\in \mathbb{F}$, the above claims
have proved $D_0(\mathfrak{L})\subset \mathbb{F} (C\otimes C)$.
Because $\mathfrak{L}=[\mathfrak{L},\mathfrak{L}]$, we obtain
$D_0(\mathfrak{L})\subset \mathfrak{L} \cdot D_0(\mathfrak{L})=0.$

\begin{clai}\adddot \label{clai5} \rm
For every $D\in{\rm Der}(\mathfrak{L},V)$, (\ref{summable}) is a
finite sum.
\end{clai}
\vskip4pt
\par According to the above results, for any $\alpha\in \Gamma,$
we can suppose $D_\alpha=(u_\alpha)_{inn}+\lambda\otimes C+C\otimes
\eta$ for some $u_\alpha\in V_\alpha,$ $\lambda,\eta\in \mathbb{F}.$
If $\{\alpha \mid u_\alpha \neq 0\}$ is an infinite set, then, we
obtain $D(L_0)= \underset{\alpha\in \Gamma}{\sum}L_0\cdot
u_\alpha=\underset{\alpha\in \Gamma}{\sum} \alpha
u_\alpha+\lambda\otimes C+C \otimes \eta$ is an infinite sum, a
contradiction with the fact that $D$ is a derivation from
$\mathfrak{L}\rightarrow V$. This proves the claim and Proposition
2.4.

Now, we can complete the proof of Theorem \ref{main} {\rm(1)} as
follows. First we need

\begin{lemm}\adddot
\label{lemma3??} Suppose $r\in V ( {\rm mod } \mathbb{F}(\mathcal
{C}\otimes \mathcal {C}))$ such that $\omega\cdot r\in {\rm
Im}(1-\tau)$ for all $\omega\in \mathfrak{L} .$ Then $r\in {\rm
Im}(1-\tau).$
\end{lemm}

\ni{\it Proof.~} First note that$\mathfrak{L} \cdot {\rm
Im}(1-\tau)\subset {\rm Im}(1-\tau).$
 Write $r=\sum_{x\in \Gamma}r_\alpha$ with $r_\alpha \in V_\alpha.$ Obviously,
 \begin{equation}
 \label{equa2.36}
 r\in {\rm Im}(1-\tau)\ \ \Longleftrightarrow \ \ r_\alpha \in {\rm Im}(1-\tau)
 \mbox{ for all } \alpha\in \Gamma.
 \end{equation}
Thus without loss of generality, we can suppose $r=r_\alpha$ is
 homogeneous.

For any $\alpha\ne0$, then $r_\alpha=\frac{1}{\alpha}L_0\cdot
r_\alpha\in {\rm Im}(1-\tau).$ Thus assume $\alpha=0.$ Now we can
write
\begin{eqnarray}\label{equa2.37}
\begin{split}
r_0=&\underset{\mu\in \Gamma/\mathbb{Z},i\in
\mathbb{Z}}{\sum}{a_{\mu,i} L_{\mu+i} \otimes L_{-\mu-i}}+
\underset{\mu\in \Gamma/\mathbb{Z},i\in \mathbb{Z}}{\sum}{b_{\mu,i}
L_{\mu+i} \otimes I_{-\mu-i}}\\
&+\underset{\mu\in \Gamma/\mathbb{Z},i\in
\mathbb{Z}}{\sum}{c_{\mu,i} I_{\mu+i} \otimes L_{-\mu-i}}+
\underset{\mu\in \Gamma/\mathbb{Z},i\in
\mathbb{Z}}{\sum}{d_{\mu,i} I_{\mu+i} \otimes I_{-\mu-i}}\\
&+ a_{\mu}L_{0}\otimes C+b_{\mu} C \otimes L_{0}+c_{\mu}I_{0}\otimes
C+d_{\mu}C \otimes I_{0},
\end{split}
\end{eqnarray}
for some
 $a_{\mu,i},b_{\mu,i},c_{\mu,i},d_{\mu,i} \in
  \mathbb{F}$, where $\{(\mu,i)\mid a_{\mu,i}\neq 0 \}$,  $\{(\mu,i)\mid b_{\mu,i}\neq 0
  \}$, $\{(\mu,i)\mid c_{\mu,i}\neq 0 \}$, $\{(\mu,i)\mid d_{\mu,i}\neq 0
  \}$ are finite sets. For a fixed $\mu \in \Gamma/\mathbb{Z},$ let
  $H_1=max\{|i| \mid a_{\mu,i}\neq 0\},H_2=min\{|i |\mid a_{\mu,i}\neq
  0\}.$  We call $L_{\mu+H_1}\otimes L_{-\mu-H_1}$ and $L_{\mu+H_2}\otimes L_{-\mu-H_2}$
 respectively the highest and the lowest terms of $r_0.$

 Choose $k\in \mathbb{Z},$ $k\neq H_1,-H_2.$ Then, $L_{k+\mu+H_1}\otimes
 L_{-\mu-H_1}$  and $L_{\mu+H_2}\otimes L_{-\mu-H_2+k}$ are the
 highest and the lowest terms of $L_k\cdot r_0$ with coefficients $(\mu+H_1-k)a_{\mu,H_1}$
and $(-\mu-H_2-k)a_{\mu,H_2}$ respectively. Since $L_k\cdot r_0\in
\mathrm{Im}(1-\tau),$ we must have $H_1=-H_2$ and
$a_{\mu,H_1}=-b_{\mu,H_2}.$ Induction on $H_1$ shows
$a_{\mu,i}=-a_{\mu,-i}$ for all $i\in\mathbb{Z}$. Similarly, we can
obtain
$b_{\mu,i}=-b_{\mu,-i}$,$c_{\mu,i}=-c_{\mu,-i}$,$d_{\mu,i}=-d_{\mu,-i}$
for all $i\in \mathbb{Z},$ and $a_\mu=-b_\mu,c_\mu=-d_\mu.$ Thus
$r_0\in \mathrm{Im}(1-\tau).$ This proves the lemma.

 \ni{\it Proof of Theorem \ref{main}(1)}
Let ($\mathfrak{L} ,[\cdot,\cdot],\D$) be a Lie bialgebra structure
on $\mathfrak{L}$. By  (\ref{deriv}), Definition \ref{def1} and
Proposition(\ref{lemma3}),$\Delta=\Delta_r+\sigma$, where $r\in
V(\mathrm{mod} (\mathcal {C}\otimes \mathcal {C}))$ and $\sigma \in
\mathbb{F}\otimes \mathbb{F}C+\mathbb{F}C \otimes \mathbb{F}$. By
(\ref{cLie-s-s}) $\mathrm{Im}\Delta \subset \mathrm{Im}(1-\tau),$ so
$\Delta_r(L_\alpha)+\sigma (L_\alpha)\in \mathrm{Im}(1-\tau)$ for
$\alpha \in \Gamma$, which implies that $c_\mu=-d_\mu$.  Moreover,
$\Delta_r(I_\alpha) \in \mathrm{Im}(1-\tau)$ for  $\alpha\in
\Gamma$. Thus, $\sigma (\mathfrak{L})\in \mathrm{Im}(1-\tau)$. So
$\mathrm{Im} \Delta_r \in \mathrm{Im}(1-\tau)$. It follows
immediately from Lemma\ref{lemma3??} that $r\in
\mathrm{Im}(1-\tau)(\mathrm{mod} (\mathcal {C}\otimes \mathcal
{C}))$. According to the above results,  we always have $0\neq
\lambda \in \mathbb{F}$ make $\sigma =\lambda \otimes C-C\otimes
\lambda.$ Then for any $\omega_\alpha \in \mathfrak{L}_\alpha,\alpha
\in \Gamma,$ we have
$$\begin{array}{l}
(1+\xi+\xi^2)\cdot (1\otimes \sigma)\cdot \sigma(\omega_\alpha)   \\[12pt]
 =\lambda(1-\delta_{\alpha,0})w_1(1+\xi+\xi^2)(1\otimes \sigma)(I_\alpha \otimes C-C\otimes
 I_\alpha)   \\[12pt]
 =0.
 \end{array}$$
It shows ($\mathfrak{L} ,[\cdot,\cdot],\sigma$) is a Lie bialgebra,
and the proof of Theorem \ref{main}(1) is completed. But by
Proposition (\ref{lemma3}), there is no $0\neq r\in \mathfrak{L},$
such that $\Delta=x\cdot r.$ Thus, ($\mathfrak{L}
,[\cdot,\cdot],\D)$ can not be a coboundary
 Lie bialgebra.

\section{Lie bialgebras of the centerless generlized
Heisenberg-Virasoro algebra}

\setcounter{section}{3} \setcounter{theo}{0}\setcounter{equation}{0}

Denote $\overline{\mathfrak{L}}=\mathfrak{L}/\mathcal {C}$, then
$\overline{\mathfrak{L}}$ is the centerless generalized
Heisenberg-Virasoro algebra.

\begin{theo}
 \label{main2} {\rm(1)} Every Lie bialgebra structure on the Lie
 algebra $\overline{\mathfrak{L}} $ is a triangular coboundary Lie
 bialgebra.

 {\rm(2)} An element $r \in \Im(1 - \tau) \subset
 {\overline{\mathfrak{L}}}\otimes \overline{{\mathfrak{L}}} $
  satisfies CYBE if and only if it
 satisfies MYBE.

 {\rm(3)} Regarding $\overline{V}=\overline{{\mathfrak{L}} }\otimes \overline{{\mathfrak{L}}}$ as an
 $\overline{{\mathfrak{L}}}$-module under the adjoint diagonal action of
 $\overline{{\mathfrak{L}}} $, we have $H^1(\overline{{\mathfrak{L}} }
 ,\overline{V})=\Der(\overline{{\mathfrak{L}}} ,\overline{V})/\Inn(\overline{{\mathfrak{L}}} ,\overline{V})=0.$
\end{theo}
  \vskip10pt

 \ni{\it Proof}\,\,  According to the  Theorem
 \ref{main}{\rm(2)}{\rm(3)}, $\Der(\overline{{\mathfrak{L}}}
,\overline{V})/\Inn(\overline{{\mathfrak{L}}} ,\overline{V})=0$. So
{\rm(2)},{\rm(3)} hold obviously. By Theorem \ref{main}{\rm(1)}, we
obtain $\Delta=\Delta_r, r\in \mathrm{Im}(1-\tau)$. Then by Lemma
\ref{some}, we have $c(r)=0$.  Thus,
$(\overline{\mathfrak{L}},[\cdot,\cdot],\Delta)$ is a triangular
coboundary Lie bialgebra. Theorem is proved.\QED

\bibliographystyle{amsplain}

\end{document}